[1996/12/01]
\documentclass{article}
\usepackage{amsmath,amsthm}
\newtheorem{thm}{Theorem}[section]
\newtheorem{lem}[thm]{Lemma}

\theoremstyle{definition}

\let\abs=\envert
\newcommand{\floor}[1]{\left\lfloor#1\right\rfloor}

\theoremstyle{remark}

\title{On the divisibility of odd perfect numbers by a high power of a prime
\footnote{The original version submitted in Apr 4, 2005.  The revised version submitted in Aug 17, 2005.}}

\date{}
\author{Tomohiro Yamada}
\begin{document}
\maketitle

\begin{abstract}
We study some divisibility properties of multiperfect numbers.  Our main result is:
if $N=p_1^{\alpha_1}\cdots p_s^{\alpha_s} q_1^{2\beta_1}\cdots q_t^{2\beta_t}$
with $\beta_1, \cdots, \beta_t$ in some finite set $S$
satisfies $\sigma(N)=\frac{n}{d}N$, then $N$ has a prime factor smaller than $C$,
where $C$ is an effective computable constant depending only on $s, n, S$.
\end{abstract}

\section{Introduction}\label{intro}
We denote by $\sigma(N)$ the sum of divisors of $N$ a positive integer and define
$h(N)=\sigma(N)/N$.  $N$ is said to be $n/d$-perfect if $h(N)=n/d$ and said to be
perfect if $h(N)=2$.  It has been known an odd perfect number must satisfy
various conditions.  Suppose $N$ is an odd perfect number.  Euler has shown that
$N=p^{\alpha} q_1^{2\beta_1}\cdots q_t^{2\beta_t}$ for distinct
odd primes $p, q_1, \cdots, q_t$ with $p\equiv \alpha\equiv 1\pmod{4}$.
Steuerwald\cite{St} proved that we cannot have $\beta_1=\cdots=\beta_t=1$.
McDaniel\cite{Mc1} proved that we cannot have $\beta_1\equiv\cdots\equiv\beta_t\equiv1\pmod{3}$.
If $\beta_1=\cdots=\beta_t=\beta$, then it is known that
$\beta\ne 2$(Kanold\cite{Ka1}), $\beta\ne 3$(Hagis and McDaniel\cite{HMD}),
$\beta\ne 5, 12, 24, 17, 62$(McDaniel and Hagis\cite{MDH}), $\beta\ne 6, 8, 11, 14, 18$(Cohen and Williams\cite{CW}).
In their paper \cite{HMD}, Hagis and McDaniel conjecture that
$\beta_1=\cdots=\beta_t=\beta$ does not occur.  The author\cite{Ymd}
proved that there are only finitely many counterexamples
for any given $\beta$.

However, if we relax the condition $\beta_1=\cdots=\beta_t=\beta$,
then the situation becomes quite different.  The simplest problem in
this direction would be whether there exists an odd perfect number of
the form $p^{\alpha}q_1^{2\beta_1}q_2^{2\beta_2}\cdots q_t^{2\beta_t}$
with $p\equiv \alpha\equiv 1\pmod{4}$ and $\beta_i\le 2$. This problem
has been studied by McDaniel\cite{Mc2} and Cohen\cite{Co}.
These papers give {\it lower} bounds for the smallest prime factor of $N$:
the former paper shows it must be $\ge 101$, and the latter shows it must be $\ge 739$.
This special case will be approached in Theorem \ref{thm2}.

In general, we can make a conjecture that for an fixed finite
set $S$ of integers, a fixed rational $n/d$ and a fixed integer s,
there exists only finitely many odd $n/d$-perfect numbers
$N=p_1^{\alpha_1}\cdots p_s^{\alpha_s} q_1^{2\beta_1}\cdots q_t^{2\beta_t}$ with $\beta_1, \cdots, \beta_t$ contained in $S$.

This conjecture still seems to be far beyond reach,
though this conjecture is weaker than the finiteness(or non-existence)
conjecture of odd $n/d$-perfect number.
In this paper, we shall show that such an odd multiperfect number, if it exists,
must have a small prime factor and we can compute an upper bound for
this prime factor in terms of $s, n, S$.

\begin{thm}\label{thm1}
Let $n, d, \beta_1, \cdots, \beta_t$ be positive integers such that $\beta_1, \cdots, \beta_t$
belong to some finite set $S$.  If $N=p_1^{\alpha_1}\cdots p_s^{\alpha_s} q_1^{2\beta_1}\cdots q_t^{2\beta_t}$
satisfies $h(N)=\frac{n}{d}$, then $N$ has a prime factor smaller than $C$,
where $C$ is an effective computable constant depending only on $s, n, S$.
\end{thm}

We use sieve method to show that $N$ has a prime factor dividing an integer determined by $S$
or the set of $q_i$'s must be thin.  In either case, we conclude that
$N$ has a prime factor smaller than $C$. The computation of $C$
requires the prime number theorem for arithmetic progression with
an effectively computable error term.

Since $C$ is effectively computable, we would be able to show there is no odd $n/d$-multiperfect number
of this form by showing any prime $<C$ could not be a divisor of an odd $n/d$-multiperfect number.  However,
there seems to be no method which is assured to determine whether a given prime can be a divisor of an odd
$n/d$-multiperfect number.  Moreover, $C$ turns out to be very large if we estimate $C$
along our method, even in (relatively) good cases for us.

We shall give an upper bound result for the smallest prime factor of $N$ by explicitly estimating $C$
in the above-mentioned special case $n/d=2$ and $S=\{2, 4\}$ in Theorem \ref{thm1}.
\begin{thm}\label{thm2}
If $N=p^e q_1^2\cdots q_s^2 q_{s+1}^4\cdots q_{s+t}^4$ is an odd perfect number,
then $N$ has a prime factor less than $\exp({4.97401\times 10^{10}})$.
\end{thm}

This upper bound is undoubtfully large, though we can make good use of the peculiarity of
the case in many steps of the proof.  Calculations of zeros and zero-free regions of
Dirichlet L-functions would improve our upper bound.  But a considerable improvement
cannot be expected.

\section{Preliminaries to Theorem \ref{thm1}}\label{lemma}
In this section, we denote by $N$ an arbitrary
positive integer.  We begin with a result concerning
the approximation of rationals by numbers of the form $h(N)$
which is interesting in itself.
This result generalizes results of Kishore\cite{Ks1}\cite{Ks2}
and is similar to a result of Pomerance \cite{Pom}.
We do not claim that this result is new, though we can
find no result of this kind in the literature.
\begin{lem}\label{lm1}
If $N=p_1^{e_1}\cdots p_k^{e_k}$ with $2<p_1<\cdots<p_k$ and $h(N)=n/d$,
then for any $s<k$ there exists an effectively computable constant $\delta>0$ depending
only on $n, s$ for which
\begin{equation}\label{eq10}
h(p_1^{e_1}\cdots p_s^{e_s})\le \frac{n}{d}-\delta
\end{equation}
holds.
\end{lem}
\begin{proof}
We begin with the case $s=1$.
First we note that $d$ must be odd since $N$ is odd.
So we have $p_1/(p_1-1)\ne n/d$.
If $p_1/(p_1-1)<n/d$, then $p_1>n/(n-d)$ and therefore
$p_1\ge (n+1)/(n-d)$.  Hence we obtain
\begin{equation}\label{eq11}
h(p_1^{e_1})\le \frac{p_1}{p_1-1}\le\frac{n+1}{d+1}\le\frac{n}{d}-\frac{1}{d(d+1)}\le\frac{n}{d}-\frac{1}{n^2}.
\end{equation}
On the other hand, if $p_1/(p_1-1)>n/d$, then we have
\begin{equation*}
\frac{p_1}{p_1-1}(1-p_1^{-e_1-1})\le h(p_1^{e_1})\le\frac{n}{d}\le \frac{p_1}{p_1-1}-\frac{1}{d(p_1-1)}=\frac{p_1}{p_1-1}(1-\frac{1}{dp_1}).
\end{equation*}
and therefore $p_1^{e_1}\le d$.  Hence we obtain
$h(p_1^{e_1})\le n/d-1/d^2$.  Combining this result with (\ref{eq11}),
we conclude (\ref{eq10}) holds for $s=1$ with $\delta=1/n^2$.

We suppose that (\ref{eq10}) holds for $s-1$ in place of $s$
with some $\delta'$ in place of $\delta$.  Then we shall show that
(\ref{eq10}) also holds for $s$ and with some $\delta$.

We may assume without loss of generality that
\begin{equation}\label{eq12}
h(p_1^{e_1}\cdots p_s^{e_s})\ge n/d-\delta'/2.
\end{equation}
Hence we have $h(p_s^{e_s})\ge (n/d-\delta'/2)/(n/d-\delta')>1$,
which implies that $p_s$ is bounded by some effectively computable
constant $C$ depending only on $n$ and $s$. Hence so are $p_1, \cdots, p_{s-1}$.
The argument of Pomerance\cite{Pom} implies that there is
an effectively computable constant $\delta(s, n, p_1, \cdots, p_s)$
if $h(N)=n/d$ and $n/d-h(p_1^{e_1}\cdots p_s^{e_s})>0$, then
$n/d-h(p_1^{e_1}\cdots p_s^{e_s})>\delta$.  Indeed, (4.5) of \cite{Pom} states that
\begin{equation}
\delta(1, n, p_1)\ge \min\{\frac{1}{n^2(p-1)}\}
\end{equation}
and page 200 of \cite{Pom} shows that we have either
\begin{equation}
\delta(s, n, p_1, \cdots, p_s)\ge\frac{1}{n\prod(p_i-1)}
\end{equation}
or
\begin{equation}
\begin{split}
\delta(s, n, p_1, \cdots, p_s)\ge&\min\{h(p_i^{a_i})\delta(s-1, n, p_1, \cdots, \hat p_i, \cdots, p_s)\\&\vert 1\le i\le s, 1\le a_i\le x_i \},
\end{split}
\end{equation}
where $x_i=\floor{\log(2sn\prod(p_i-1))/\log p_i}$.

\end{proof}

We need the prime number theorem for arithmetic progressions with an effectively computable error term.

\begin{lem}\label{lm2}
Let $\pi(N, d, a)$ be a number of primes up to $N$ which is congruent to $a\pmod{d}$.  If $(d, a)=1$, then
\begin{equation}
\pi(N, d, a)=\frac{N}{\varphi(d)\log{N}}+O(\frac{N}{\log^2{N}}),
\end{equation}
where the implied constant is effectively computable in terms of $d$ and $a$.
\end{lem}

This lemma follows from Theorem 9.6 in Karatsuba \cite{Krt}.  Another result that we need is a standard result in sieve theory.

\begin{lem}\label{lm3}
Let $A$ and $\Omega_p$, where $p$ is an arbitrary prime number, be sets of positive integers,
$B$ be a positive integer, $X$ be a real number, and $\rho(n)$ be a multiplicative arithmetic function satisfying
$0\le\rho(p)\le\min\{B, p-1\}$ for any prime $p$.  Denote by $A_d$ the set of positive integers in $A$
which belongs to $\Omega_p$ for any $p$ dividing $d$.  Define
\begin{equation}
P(z)=\prod_{p<z, p\text{ is prime}}p.
\end{equation}
\begin{equation}
R_d=\frac{\rho(d)}{d}X-\#A_d
\end{equation}
and
\begin{equation}
S(A, P(z))=\#\{a\in A: a\text{ does not belong to }\Omega_p\text{ for any prime }p\mid P(z)\}.
\end{equation}
If $0<s<1/2$, $z=X^s$ and $\abs{R_d}<d$ for $d<z$, then
\begin{equation}
S(A, P(z))=O(X\prod_{p<z}(1-\frac{\rho(p)}{p})),
\end{equation}
where the implied constant is effectively computable in terms of $B$ and $s$.
\end{lem}

This lemma follows from the Brun-Selberg Sieve method.
This lemma is a generalizaion of Corollary 2.2.1.1 or Corollary 3.3.1.2 in \cite{Gre}.  We can easily extend
these results to the theorem mentioned above(See section 1.3.4 in \cite{Gre}).

There are several explicit upper bound sieve formula to obtain explicit upper bound for the implied constant in this Lemma.
We use the upper bound formula given in Theorem 2.1.1 and Theorem 2.2.1 in \cite{Gre}.
We begin with defining
\begin{equation*}
g(p)=\frac{\rho(p)}{p-\rho(p)},
\end{equation*}
\begin{equation*}
E(D, P)=\sum_{d_1, d_2<\sqrt{D}, d_1, d_2\mid P}\abs{R_{[d_1, d_2]}},
\end{equation*}
\begin{equation*}
V(P(z))=\prod_{p<z, p\text{ is prime}}(1-\frac{\rho(p)}{p})
\end{equation*}
\begin{equation*}
G_w(x)=\sum_{d\le x, d|P(w)}g(d).
\end{equation*}
The following three lemmas concern the upper bound sieve inequality.
These inequalities allows us to calculate an upper bound in Theorem \ref{thm1} explicitly.
\begin{lem}\label{lm21}
Let $z\le w\le \sqrt{D}$.  Then
\begin{equation}
S(A, P(w))\le \frac{X}{G_w(\sqrt{D})}+E(D, P(w)).
\end{equation}
\end{lem}
\begin{proof}
This follows from Theorem 2.1.1, Corollary 2.1.2.1, and (2.2.1.3) in \cite{Gre}.
\end{proof}
\begin{lem}\label{lm22}
If $R_d\le\rho(d)$ for each $d$ dividing $P(w)$, then
\begin{equation}
E(D, P(w))\le D\prod_{p\mid P(w)}(1+\frac{\rho(p)}{p})^3.
\end{equation}
\end{lem}
\begin{proof}
We observe that each $d$ dividing $P(w)$ has at most $3^{\omega(d)}$ representations in the form $[d_1, d_2]$.
We shall omit the rest of the proof since it proceeds as in pages 100-101 in \cite{HR}.
\end{proof}
\begin{lem}\label{lm23}
If
\begin{equation}
\frac{1}{\log{z}}\sum_{p<z}\frac{\rho(p)\log{p}}{p}\le B
\end{equation}
holds, then
\begin{equation}
G_w(\sqrt{D})\ge \psi_0(v)/V(P(z)),
\end{equation}
where
\begin{equation}
\psi_0(v)=1-\exp(-\psi(B, v/2))
\end{equation}
with $\psi(B, v)$ defined by
\begin{equation}
\psi(B, v)=\max\{0, v\log{\frac{v}{B}}-v+B\}.
\end{equation}
\end{lem}
\begin{proof}
This is Theorem 2.2.1 in \cite{Gre}
\end{proof}

\section{Proof of Theorem \ref{thm1}}
Let $N=p_1^{\alpha_1}\cdots p_s^{\alpha_s} q_1^{2\beta_1}\cdots q_t^{2\beta_t}$ be a solution of $h(N)=\frac{n}{d}$.
Let us denote by $P$ the set $\{p: p\mid(2\beta+1)\text{ for some }\beta\in S\}$.  Then $P$ is a finite set of primes.
Let us denote by $P$ also the product $\prod_{p\in P}p$ and let $T=\{i: q_i\equiv 1\pmod{P}\}$.
We define $\Omega_P(N)$ to be the number of prime factors of $N$ which belongs to $P$, counting multiplicity.
Since $P$ is a finite set depending only on $S$, we may assume without loss of generality
that $N$ has no prime divisor in $P$ so that $\Omega_P(N)=0$.
We denote $Q_{\beta}$ by the set of primes $q_i$ with $\beta_i=\beta$.

We shall begin the proof with an simple observation.  There are at most $\Omega_P(nN)=\Omega_P(n)$ prime factors $q_i$ in $T$ since if $q_i\in T$, then $\sigma(q_i^{2\beta_i})$ is divisible by $2\beta_i$.  Therefore the number of prime factors of $\sigma(N)=nN/d$ congruent to $1\pmod{P}$ is at most $s+\Omega_P(n)+\omega(n)$.  Denote by $U$ the set of primes $\equiv 1\pmod{P}$ not dividing $\sigma(N)$.  Hence we see that $U$ is a set of primes $\equiv 1\pmod{P}$ in $T$ except at most $s+\Omega_P(n)+\omega(n)$ primes.  This allows us to apply the sieve method described in the previous section.

Now we prove a result concerning the distribution of prime factors of $N$, which is the most important lemma
in the proof of Theorem \ref{thm1}.
\begin{lem}\label{lm32}
There exist effectively computable constants $c, \kappa$ depending only $s$ and $n$ such that
\begin{equation}
\#\{q: q\in Q_{\beta}, q\le x\}\le c\frac{x}{(\log{x})^{1+\kappa}}.
\end{equation}
\end{lem}
\begin{proof}
We denote by $B_0, B_1, \cdots$ effectively computable constants depending only on $P$.
If $q\in Q_{\beta}$, then $q$ is prime and $\sigma(q_i^{2\beta})$ has no prime factor in $U$.
Let $p$ be a prime factor of ${2\beta+1}$.  Then $p\mid P$.  Thus, if $r\in U$, then $r$ is
a prime congruent to $1\pmod{p}$.  Hence there are $p-1$ congruent classes $g_1(r), \cdots, g_{p-1}(r) \pmod{r}$
belonging to order $p$.
Clearly, $r$ does not divide $\sigma(q_i^{p-1})$.  Hence $q_i$ belongs to none of $p$ classes $0, g_1, \cdots, g_{p-1} \pmod{r}$.

Now we can apply the sieve method described in the previous section with $A$ the set of integers $\le x$, $X=x$,
$\Omega_r$ the set of integers $\le x$ belongs to any of congruent classes $0, g_1, \cdots, g_{p-1}\pmod{r}$
for $r\in U$ and $0\pmod{r}$ for $r\not\in U$, $\rho(r)=p$ for $r\in U$ and $\rho(r)=1$ for $r\not\in U$.

By Lemma $\ref{lm2}$, we have for $z>B_0$,
\begin{equation}\label{eq31}
\sum_{r\le z}\frac{\rho(r)\log{r}}{r}\le\sum_{r\le z}\frac{\log{r}}{r}+\sum_{r\le z, r\in U}\frac{(p-1)\log{r}}{r}\le(B_1+(p-1)B_2)\log{z},
\end{equation}
\begin{equation}\label{eq311}
\sum_{r\le z, r\equiv 1\pmod{P}}\frac{p-1}{r}\ge(p-1)(\frac{\log{\log{z}}}{\phi(P)}-B_3).
\end{equation}
\begin{equation}\label{eq312}
\sum_{r\le z, r\in U}\frac{p-1}{r}\le(p-1)(\frac{\log{\log{z}}}{\phi(P)}+B_4),
\end{equation}

We use the well-known formula of Mertens and recall that
the number of primes $r\equiv 1\pmod{P}$ not contained in $U$
is finite and explicitly computable to obtain
\begin{equation}\label{eq32}
\begin{split}
V(P(z))=\prod_{r<z}(1-\frac{\rho(r)}{r})&\le\prod_{r<z}(1-\frac{1}{r})\prod_{r<z, r\in U}(1-\frac{1}{r})^{p-1}\\
&\le\prod_{r<z}(1-\frac{1}{r})\prod_{r<z, r\equiv 1\pmod{P}}(1-\frac{1}{r})^{p-1}\\
&\times\prod_{r\equiv 1\pmod{P}, r\not\in U}(1-\frac{1}{r})^{-(p-1)}\\
&\le (e^{B_3}B_5)^{p-1}B_6(\log{z})^{-(1+\kappa)},
\end{split}
\end{equation}
\begin{equation}\label{eq33}
\prod_{r<z}(1+\frac{\rho(r)}{r})\le\prod_{r<z}(1+\frac{1}{r})\prod_{r<z, r\in U}(1+\frac{1}{r})^{p-1}\le e^{B_4(p-1)}B_7(\log{z})^{1+\kappa},
\end{equation}
where $\kappa=\frac{p-1}{\phi(P)}$, since $-\log{(1-z)}\ge z$ for $0\ge z<1$ and $1+w\le e^{w}$ for $w\le 0$.
We note that $B_0, \cdots, B_7$ depend only on $P$.

We put $D=X/(\log{X})^{4(1+\kappa)}$ and  $z=w=D^{1/v}$.  We denote by $\psi_1(v), \cdots$ real-valued functions of $v$
depending only on $P$, $p$ and $D$.
If $X\ge \psi_3(v)$, then $z>B_0$ and (\ref{eq32}) yields
\begin{equation}\label{eq34}
V(P(z))\le \psi_1(v)(\log{z})^{-(1+\kappa)}\le \psi_2(v)(\log{X})^{-(1+\kappa)}.
\end{equation}

Hence, by Lemma \ref{lm23}, we have
\begin{equation}\label{eq35}
\frac{1}{G_w(\sqrt{D})}\le \frac{1}{G_z(\sqrt{D})}\le \frac{\psi_2(v)}{\psi_0(v)}(\log{X})^{-(1+\kappa)}.
\end{equation}
Furthermore, by virtue of Lemma \ref{lm22}, if $X\ge \psi_3(v)$, then (\ref{eq33}) yields
\begin{equation}\label{eq36}
E(D, P(w))\le \psi_4(v)D(\log{w})^{3(1+\kappa)}\le \psi_4(v)D(\log{D})^{3(1+\kappa)}\le \psi_4(v)\frac{X}{(\log{X})^{1+\kappa}}.
\end{equation}

By Lemma \ref{lm21}, we obtain
\begin{equation}
S(A, P(z))\le(\frac{\psi_2(v)}{\psi_0(v)}+\psi_4(v))\frac{X}{(\log{X})^{1+\kappa}}=\psi_5(v)\frac{X}{(\log{X})^{1+\kappa}}.
\end{equation}

The lemma easily follows noting that $\#\{q: q\in Q_{\beta}, q\le x\}\le S(A, P(\sqrt{x}))+\sqrt{x}$.
\end{proof}

Now we shall prove Theorem \ref{thm1}.  Since
\begin{equation*}
\prod_{i=1}^{s}h(p_i^{\alpha_i})<n/d,
\end{equation*}
we obtain
\begin{equation}\label{eq13}
\prod_{i=1}^{s}h(p_i^{\alpha_i})\le n/d-\delta
\end{equation}
by Lemma \ref{lm1}.
where $\delta$ is an effectively computable constant depending on $n, S$
Let $d_{\beta}=\prod_{q} q/(q-1)$, where $q$ runs all primes in $Q_{\beta}$.  It follows from (\ref{eq13}) that $\prod_{\beta\in S}d_{\beta}$ must be $\ge \frac{n}{n-d\delta}$.  Hence we have that some $d_{\beta}\ge\mu$, where $\mu>1$ is effectively computable in terms of $s, n, S$.
Now it immidiately follows from Lemma \ref{lm32} that $Q_{\beta}$ has an element smaller than $C$.  This proves the theorem.
\section{Preliminaries to Theorem \ref{thm2}}

For our purpose, it suffices to calculate $B_0, B_1, \cdots$ in the previous section.
This requires many inequalities involving sums or products of primes in 
arithmetic prograsions.

We set $R=964.5908801$ and take $B_0$ as $e^R$. We begin by calculating $B_5$.
Noting that $n/d=2$ and $P=\{3, 5\}$,
we obtain $\Omega_P(n)=0$ and therefore $T=\phi$.  Moreover no prime
$\equiv 1\pmod{P}$ divides $n$.  Hence $U$ contains all prime $\equiv 1\pmod{P}$
with at most one exception.  Thus we can take $B_5=\frac{31}{30}$.

We refer some inequalities involving primes.
\begin{lem}\label{lm41}
For any $z>1$ we have
\begin{equation}\label{eq411}
\sum_{p\le z}\frac{\log{p}}{p}\le \log{z},
\end{equation}
\begin{equation}\label{eq412}
\prod_{p\le z}(1-\frac{1}{p})\le \frac{e^{-\gamma}}{\log{z}}(1+\frac{1}{2\log^2{z}}).
\end{equation}
Moreover, for any $z>286$ we have
\begin{equation}\label{eq413}
\prod_{p\le z}(1+\frac{1}{p})\le \frac{e^K}{\log{z}}(1+\frac{1}{\log^2{z}}),
\end{equation}
where $K$ is a constant less than $0.2615$.
\end{lem}
\begin{proof}
The inequality (\ref{eq411}) is the formula (3.24) in \cite{RS} and the inequality (\ref{eq412}) is Theorem 7 in \cite{RS}.
The inequality (\ref{eq413}) follows from (3.18) in \cite{RS}.
\end{proof}
This lemma allows one to take $B_1=1$, $B_6=e^{-\gamma}(1+(2\log^2{z})^{-1})$ and $B_7=e^K(1+(\log{z})^{-2})$.

As can be seen by the proof of Theorem \ref{thm1}, we need some results on the distribution of
prime numbers $\equiv 1\pmod{15}$ in order to estimate the constant in problem.  The starting point
is the following result due to Ramar{e} and Rumely\cite{RR} and Dusart\cite{Du}.
\begin{lem}\label{lm42}
\begin{equation}
\abs{\theta(x, 15, 1)-\frac{x}{8}}\le\frac{0.609x}{\log{x}}
\end{equation}
for any positive $x$, and
\begin{equation}
\abs{\theta(x, 15, 1)-\frac{x}{8}}\le\frac{0.008634}{8}x
\end{equation}
for any $x\ge 10^{10}$.
\end{lem}
\begin{proof}
Assume first that $x\ge e^R$.  We use Theorem 5 in \cite{Du}.  We obtain the following estimates for $X_i$ in this theorem as follows:
$8.33<X_0<8.34, 6.20<X_1<6.21, 2.72<X_2<2.73, 3.50<X_3<3.51, X_4=10$.  Now Theorem 5 in \cite{Du} gives
\begin{equation}
\abs{\theta(x, 15, 1)-\frac{x}{8}}\le\frac{0.189x}{\log{x}}\le\frac{0.008634}{8}x,
\end{equation}
for $x$ in the assumed range.

Assume next that $10^{10}\le x\le e^R$.  In three cases $x\ge 10^{100}$, $10^{30}\le x\le 10^{100}$, $x\le 10^{30}$,
Theorem 1 and Table 1 in \cite{RR} shows that the absolute value in lemma is at most
$\frac{0.005045}{8}x$, $\frac{0.007088}{8}x$, $\frac{0.008634}{8}x$ respectively.
In any case, this does not exceed $\frac{0.609x}{\log{x}}$.

Assume last that $x\le 10^{10}$.  If $x\ge 50$, then the Table 2 in \cite{RR} gives
$\abs{\theta(x, 15, 1)-\frac{x}{8}}\le 1.098\sqrt{x}\le \frac{0.609x}{\log{x}}$.
If $x\le 50$, then $\abs{\theta(x, 15, 1)-\frac{x}{8}}\le x/8\le\frac{0.609x}{\log{x}}$.
\end{proof}
\begin{lem}\label{lm43}
For any $z>e^R$ we have
\begin{equation}
\frac{1}{8}\log\log{z}-0.312\le\sum_{p\le z, p\equiv 1\pmod {15}}\frac{1}{p}\le \frac{1}{8}\log\log{z}+0.0572.
\end{equation}
Thus we can take $B_3=0.312$ and $B_4=0.0572$.
\end{lem}
\begin{proof}
Put $\phi_1(z)=\int_{60}^z\frac{\theta(t, 15, 1)}{t^2\log^2{t}}(1+\log{t})dt$.  Then, by Lemma \ref{lm42} we obtain
\begin{equation}
\begin{split}
\phi_1(e^R)&=\int_{60}^{e^R}\frac{\theta(t, 15, 1)}{t^2\log^2{t}}(1+\log{t})dt\\
&\le\int_{60}^{e^R}\frac{1}{8}(\frac{1}{t\log^2{t}}+\frac{1}{t\log{t}})+0.609(\frac{1}{t\log^3{t}}+\frac{1}{t\log^2{t}})\\
&\le\frac{1}{8}(\log{R}-\log\log{60})+0.883(\frac{1}{\log{60}}-\frac{1}{R})
\end{split}
\end{equation}
and
\begin{equation}\label{eq431}
\begin{split}
\phi_1(z)&=\phi_1(e^R)+\int_{e^R}^z\frac{\theta(t, 15, 1)}{t^2\log^2{t}}(1+\log{t})dt\\
&\le\phi_1(e^R)+\int_{e^R}^z\frac{1}{8}(\frac{1}{t\log^2{t}}+\frac{1}{t\log{t}})+0.189(\frac{1}{t\log^3{t}}+\frac{1}{t\log^2{t}})\\
&\le\phi_1(e^R)+\frac{1}{8}(\log\log{z}-\log{R})+0.315(\frac{1}{R}-\frac{1}{\log{z}})\\
&\le\frac{1}{8}(\log\log{z}-\log\log{60})+\frac{0.883}{\log{60}}-\frac{0.568}{R}.
\end{split}
\end{equation}
The sum in the lemma can be estimated by $\phi_1(z)$ as follows:
\begin{equation}\label{eq432}
\sum_{p\le z, p\equiv 1\pmod {15}}\frac{1}{p}=\frac{\theta(z, 15, 1)}{z\log{z}}-\frac{\theta(60, 15, 1)}{60\log{60}}+\frac{1}{31}+\phi_1(z).
\end{equation}
Combining (\ref{eq431}) and (\ref{eq432}), we obtain the second inequality in the lemma. The first inequality can be obtained in a similar way.
\end{proof}
\begin{lem}\label{lm44}
For any $z>e^R$ we have
\begin{equation}
\sum_{p\le z, p\equiv 1\pmod {15}}\frac{\log{p}}{p}\le 0.12615\log{z},
\end{equation}
that is, we can take $B_2=0.12615$.
\end{lem}
\begin{proof}
Put $\phi_2(z)=\int_{31}^z\frac{\theta(t, 15, 1)}{t^2}dt$.  Then we obtain
\begin{equation}
\begin{split}\label{eq441}
\phi_2(z)=&\int_{31}^{10^{10}}+\int_{10^{10}}^{e^R}+\int_{e^R}^{z}\frac{\theta(t, 15, 1)}{t^2}dt\\
&\le\frac{1}{8}(\log\frac{10^{10}}{31})+2.196(\frac{1}{\sqrt{31}}-\frac{1}{10^5})\\
&+\frac{1.008634}{8}(R-\log{10^{10}})+(\frac{1}{8}+0.000196)(\log{z}-R)\\
&\le(\frac{1}{8}+0.000196)\log{z}-\frac{1}{8}\log{31}+2.196(\frac{1}{\sqrt{31}}-\frac{1}{10^5})\\
&+\frac{0.008634}{8}(R-\log{10^{10}})-0.000196R.
\end{split}
\end{equation}
The sum in this lemma can be estimated as follows:
\begin{equation}\label{eq442}
\sum_{p\le z, p\equiv 1\pmod {15}}\frac{\log{p}}{p}=\frac{\theta(z, 15, 1)}{z}+\phi_2(z)\le \frac{1}{8}+\frac{0.189}{\log{z}}+\phi_2(z).
\end{equation}
(\ref{eq441}) and (\ref{eq442}) give
\begin{equation}
\begin{split}
\sum_{p\le z, p\equiv 1\pmod {15}}\frac{\log{p}}{p}&\le(\frac{1}{8}+0.000196)\log{z}+0.918\\
&\le(\frac{1}{8}+0.00115)\log{z}\\
&=0.12615\log{z}.
\end{split}
\end{equation}
\end{proof}

\begin{lem}\label{lm45}
Let $g$ be an integer in $\{2, 4\}$.
Denote by $\pi_g(X)$ the number of prime factors $\le X$ of $N$ with exponent $g$.  If $X>e^{14R}$, then we have
\begin{equation}
\pi_2(X)\le 40.8778\frac{X}{(\log{X})^{5/4}},
\end{equation}
and
\begin{equation}
\pi_4(X)\le 185.976\frac{X}{(\log{X})^{3/2}}.
\end{equation}
\end{lem}
\begin{proof}
To estimate $\pi_g(X)$, we use the inequality
\begin{equation}
\pi_g(X)\le \sqrt{X}+S(A, P(w))
\end{equation}
with $\rho$ defined by
\begin{equation}
\rho(p)=1
\end{equation}
and
\begin{equation}
\rho(p)=1+g.
\end{equation}
By (\ref{eq411}) and Lemma \ref{lm44}, $B$ can be taken as $1+0.12615g$.
Hence, by virtue of (\ref{eq32}), Lemma \ref{lm41} and Lemma \ref{lm43}, we can take
\begin{equation}
\psi_1(v)=(\frac{31}{30})^g e^{0.312g-\gamma}(1+\frac{1}{2}(\frac{v}{\log{D}})^2).
\end{equation}
From (\ref{eq33}), Lemmas \ref{lm22}, \ref{lm41} and \ref{lm43}, we can take
\begin{equation}
\psi_4(v)=e^{3(K+0.0572g)}v^{-3(1+g/8)}(1+(\frac{v}{\log{D}})^2)^3.
\end{equation}
We have $D(\log{D})^{3(1+g/8)}\le X(\log{X})^{-(1+g/8)}$ and $(v/\log{D})^{1+g/8}\le\\ (2v)^{1+g/8}(\log{X})^{-(1+g/8)}$.  Hence we can take
\begin{equation}
\psi_2(v)=(2v)^{1+g/8}\psi_1(v).
\end{equation}
The Case $g=2$: $\psi_5(7.019)\le 40.9177$.
The Case $g=4$: $\psi_5(7.536)\le 187.083$.
We can obtain a trivial estimate for $\sqrt{D}$ as follows:
\begin{equation}
\sqrt{D}\le D/\sqrt{D}\le X(\log{X})^{-(1+g/8)}/\sqrt{D}\le 10^{-6}X(\log{X})^{-(1+g/8)}.
\end{equation}
Finally, we can easily confirm that $z>e^R$ under the condition $X>e^{14R}$.  This completes the proof.
\end{proof}
\begin{lem}\label{lm46}
Let $\alpha$, $c$, $x_0$, $w$ be positive real numbers, $P$ be a set of primes and denote by $\pi_P(x)$ the number of primes $\le x$ in $P$.  If
\begin{equation}\label{eqe11}
\pi_P(x)\le c\frac{x}{\log^{1+\alpha}{x}}
\end{equation}
for $x>x_0$, then
\begin{equation}\label{eqe12}
\sum_{p\in P, w\le p}\frac{1}{p}\le\frac{c}{\alpha\log^{\alpha}{w}}.
\end{equation}
for $w\le x_0$.
\end{lem}
\begin{proof}
\begin{equation}
\begin{split}
\sum_{p\in P, w\le p}\frac{1}{p}&\le\int_{w}^{\infty}{\pi_P(t)}{t^2}dt\\
&\le c\int_{w}^{\infty}\frac{dt}{t\log^{1+\alpha}{t}}\\
&\le\frac{c}{\alpha\log^{\alpha}{w}}.
\end{split}
\end{equation}
\end{proof}

\section{Proof of Theorem \ref{thm2}}
We may assume that $N$ has no prime factor less than $w$ for some $w>e^{14R}$.
$N$ can be decomposed into the form $p^a KL$, where $p\equiv a\equiv 1\pmod{4}$,
$K=K^{\prime 2}$ and $L=L^{\prime 4}$ with $K^\prime$ and $L^\prime$ squarefree.
Since $h(p)<p/(p-1)\le e^{14R}/{e^{14R}-1}$, we have either $h(K)>(2-2e^{-14R})^{1/2}$
or $h(L)>(2-2e^{-14R})^{1/2}$.  Let $P=K^\prime$ if $K$ satisfies this inequality and $P=L^\prime$ otherwise.
By Lemma \ref{lm46} we obtain
\begin{equation}
\begin{split}
\sum_{p\mid P, p\ge w}\log h(p^e)&<\sum_{p \mid P, p\ge w}\log{p/(p-1)}\le\sum_{p \mid P, p\ge w}1/(p-1)\\
&\le\sum_{p \mid P, p\ge w}1/p+1/p(p-1)\\
&\le\sum_{p \mid P, p\ge w}1/p+\sum_{p\ge w}1/p(p-1)\\
&\le\sum_{p \mid P, p\ge w}1/p+\frac{w}{w-1}\sum_{p\ge w}1/p^2\\
&\le\sum_{p \mid P, p\ge w}1/p+\frac{2}{w-1}\le\frac{c}{\alpha}{\log^{\alpha}{w}}+\frac{2}{w-1}.
\end{split}
\end{equation}
Hence we have $\frac{c}{\alpha}{\log^{\alpha}{w}}>\frac{1}{2}\log(2-2e^{-14R})-\frac{2}{w-1}>0.346573$.  Therefore we obtain
\begin{equation}
w<\exp({(c/0.346573\alpha)^{1/\alpha}})
\end{equation}
Now, by Lemma \ref{lm45}, $(\alpha, c)$ can be taken as $(1/2, 187.083)$ or $(1/4, 40.9177)$. Hence we conclude that
\begin{equation}
w<\exp({4.97401\times 10^{10}}).
\end{equation}
{}
\vskip 12pt

{\small Tomohiro Yamada}\\
{\small Department of Mathematics\\Faculty of Science\\Kyoto University\\Kyoto 606-8502\\Japan}\\
{\small e-mail: \protect\normalfont\ttfamily{tyamada@math.kyoto-u.ac.jp}}
\end{document}